\documentclass[10pt]{amsart}
\usepackage[letterpaper,margin=1in,footskip=0.25in]{geometry}

\usepackage{microtype}
\usepackage{eucal}
\usepackage{setspace}
\usepackage{mathrsfs}
\usepackage{quiver}

\usepackage[french, english]{babel}

\usepackage[utf8x]{inputenc}

\usepackage{amsmath,amssymb,amsthm}
\usepackage{mathtools}
\makeatletter
  \newcommand{\supsize}{%
    \expandafter\ifx\csname S@\f@size\endcsname\relax
      \calculate@math@sizes
    \fi
    \csname S@\f@size\endcsname
    \fontsize\sf@size\z@\selectfont
  }
  \DeclareRobustCommand{\tsup}[1]{%
    \leavevmode\raise.9ex\hbox{\supsize #1}%
  }
  \DeclareTextSymbolDefault{\textprimechar}{OMS}
  \DeclareTextSymbol{\textprimechar}{OMS}{48}
  \DeclareRobustCommand{\tprime}{\tsup{\textprimechar}}
  \ProvideTextCommandDefault{\cprime}{\tprime}
\makeatother

\usepackage{mathabx}

\usepackage[shortlabels]{enumitem}
\setlist{noitemsep}

\usepackage[pdfusetitle,colorlinks]{hyperref}
\hypersetup{allcolors=blue}
\usepackage[capitalise,noabbrev]{cleveref}

\usepackage{colonequals}

\usepackage{graphicx}
\usepackage{extarrows}
\usepackage{tikz,tikz-cd}

\crefformat{equation}{\ensuremath{(#2#1#3)}}
\crefmultiformat{equation}{\ensuremath{(#2#1#3)}}{ and~\ensuremath{(#2#1#3)}}{, \ensuremath{(#2#1#3)}}{, and~\ensuremath{(#2#1#3)}}

\numberwithin{figure}{section}
\numberwithin{equation}{section}
\theoremstyle{definition}

\newtheorem{lemma}[figure]{Lemma}

\newtheorem{example}[figure]{Example}

\newtheorem{thm}[figure]{Theorem}

\newtheorem{cor}[figure]{Corollary}

\newtheorem{defn}[figure]{Definition}

\newtheorem{quest}[figure]{Question}

\newtheoremstyle{cited}{.5\baselineskip\@plus.2\baselineskip\@minus.2\baselineskip}{.5\baselineskip\@plus.2\baselineskip\@minus.2\baselineskip}{\itshape}{}{\bfseries}{\bfseries .}{5pt plus 1pt minus 1pt}{\thmname{#1}\thmnumber{ #2}\thmnote{ \normalfont #3}}
\theoremstyle{cited}

\newcommand{\Nil}{\mathrm{Nil}}

\newcommand{\Spec}{\operatorname{Spec}}
\newcommand{\Spf}{\operatorname{Spf}}

\newcommand{\Spa}{\operatorname{Spa}}

\newcommand{\ins}{\subset}

\newcommand{\dimtr}{\mathrm{dim.tr}}
\newcommand{\Frac}{\mathrm{Frac}}
\newcommand{\perfd}{\operatorname{perfd}}
\newcommand{\spc}{\mathrm{sp}}

\providecommand{\trdeg}{\operatorname{tr.deg}}
\providecommand{\toptrdeg}{\operatorname{top.tr.deg}}
\providecommand{\Toptrdeg}{\operatorname{Top.tr.deg}}
\providecommand{\trc}{\operatorname{tr.c}}

\begin{document}

\title{Topological transcendence degree and the dimension problem for perfectoid spaces}
\author{Dimitri Dine, Jack J Garzella}
\maketitle

\begin{abstract}
    
    Let $K$ be a perfectoid field of characteristic $p$ with countable value group (for example, $K=\mathbb{C}_{p}^{\flat}$). We study Zariski-closed subsets \begin{equation*}Z\subseteq \mathbb{B}_{K}^{n,\perfd}\end{equation*}of the $n$-dimensional perfectoid closed unit ball $\mathbb{B}_{K}^{n,\perfd}$ over $K$. We use an example of Gleason, and its generalization obtained in the authors' previous work, to show that two natural notions of dimension for such Zariski-closed subsets, the Krull dimension $\dim(Z)$ and Huber's iterated transcendence degree $\dimtr(Z)$, do not coincide, and that the gap between them can be arbitrarily large. 
\end{abstract}

\section{Introduction}

Fix a prime number $p$. In \cite{Scholze, Fontaine, kedlaya-liu-2015-relative-p-adic, Kedlaya-Liu2, fargues-scholze-2021-geometrization, Diamonds}, Scholze and others develop the
theory of perfectoid spaces. One of the major questions about perfectoid rings and spaces is whether their theory includes
a non-noetherian analogue of ``finite type geometry'', in a similar way that
topologically finite type algebras and rigid-analytic varieties provide an analogue of finite type geometry.

One important part of finite type geometry would be dimension theory (for example, see
\cite[Problem~I.11.1]{fargues-scholze-2021-geometrization}). In both scheme-theoretic algebraic geometry
and rigid geometry, one has a ``dimension theory'', i.e., a few different notions of dimension and
a theorem saying that they are all equal. In the perfectoid context, an analogue of affine space in
algebraic geometry and of an $n$-dimensional closed unit ball in rigid geometry is given by the
$n$-dimensional perfectoid closed unit ball \begin{equation*}\mathbb{B}_{K}^{n,\perfd}=\Spa(T_{n,
    K}^{\perfd}, T_{n, K}^{\perfd\circ})\end{equation*}over a perfectoid field $K$, whose ``ring of
    functions'' is the perfectoid Tate algebra \begin{equation*}T_{n, K}^{\perfd}=K\langle
        X_{1}^{1/p^{\infty}},\dots, X_{n}^{1/p^{\infty}}\rangle,\end{equation*}a natural analogue of
        the usual Tate algebra of convergent power series $T_{n, K}=K\langle X_{1},\dots,
        X_{n}\rangle$. In the case of a Zariski-closed subset \begin{equation*}Z=\Spa(A,
        A^{+})\subseteq \mathbb{B}_{K}^{n,\perfd}\end{equation*}we have at least four notions which
        could be reasonable candidates for a definition of dimension:
        
\begin{enumerate}[(i)]
    \item The Krull dimension of the special fiber of a(ny) formal model (as a scheme)
    \item The Krull dimension of the adic space \(Z=\Spa(A, A^{+})\) as a spectral space, called
        \(\dim(Z)\)
    \item The iterated transcendence degree, defined by Huber and denoted \(\dimtr(Z)\)
    \item The ($\ell$-adic) cohomological dimension.
\end{enumerate}

In this paper, we show that (ii) and (iii) do not coincide, and indeed can be arbitrarily far apart.

\begin{thm}[Theorem~\ref{Main theorem 3}]\label{thm:dim:noteq}
Let $K$ be a perfectoid field of characteristic $p$ whose value group is countable. 
For every $n\geq1, m\geq0$, there exists a Zariski-closed subset \begin{equation*}Z\subseteq \mathbb{B}_{K}^{n+m+2,\perfd}\end{equation*}of the $(n+m+2)$-dimensional
perfectoid closed unit ball $\mathbb{B}_{K}^{n+m+2,\perfd}$ over $K$ such that $\dim(Z)=m$ and
$\dimtr(Z)\geq n+m$.\end{thm}
The proof of the theorem follows very quickly from a fundamental example due to Gleason
\cite{gleason-2022-perfectoid-Nullstellensatz} (and its generalization due to the authors in \cite{Dine-Garzella25}), which shows
that perfectoid Tate algebras do not satisfy an analogue of the Nullstellensatz.
Let \((r_{1}, \ldots, r_{n})\) be a ($K$-)free polyradius, i.e., a tuple of positive real numbers which are linearly independent as elements of the ``multiplicative'' $\mathbb{Q}$-vector space \begin{equation*}(\mathbb{R}_{>0}/\sqrt{\vert K^{\times}\vert})\otimes_{\mathbb{Z}}\mathbb{Q}.\end{equation*}
To such a polyradius, Berkovich associates a nonarchimedean field \(K_{r_{1}, \ldots, r_{n}}\),
the completion of the field \(K(x_{1}, \ldots, x_{n})\) 
under the \((r_{1}, \ldots, r_{n})\)-Gauss norm.
Let \(K_{r_{1}, \ldots, r_{n}}^{\perfd}\) denote the field obtained
by adjoining \(p\)-power roots of all variables and further completing
(equivalently, adjoining $p$-power roots of the variables in $K(x_{1},\dots, x_{n})$ before completing).
Gleason shows:
\begin{thm}[\cite{gleason-2022-perfectoid-Nullstellensatz}, Theorem~4.1 and its proof]
Let $K$ be an algebraically closed perfectoid field of characteristic $p$ such that its value group is not equal to all of $\mathbb{R}_{>0}$, 
and let $r\in (0, 1]\setminus \vert K^{\times}\vert$. Then the nonarchimedean field \(K_{r}^{\perfd}\) is a quotient of \(T_{3, K}^{\perfd}\).
\end{thm}

In \cite{Dine-Garzella25}, the authors generalized this to:
\begin{thm}
    [Theorem~\ref{Old theorem}, \cite{Dine-Garzella25}, Theorem~6.1]
    \label{thm:kr:bigquotient}Let $K$ be a perfectoid field of characteristic $p$ and let $(r_1,\dots, r_n)\in (0, 1]^{n}$ be a $K$-free polyradius. Then the nonarchimedean field
    \(K_{r_{1}, \ldots, r_{n}}^{\perfd}\) is a quotient of \(T_{n+2, K}^{\perfd}\). 
    
    Moreover, the quotient map restricts to a surjective map from the subring of power-bounded elements $T_{n+2, K}^{\perfd\circ}$ onto $K_{r_{1},\ldots, r_{n}}^{\perfd\circ}$.
\end{thm}
Now, the main theorem follows from this and (a variant of) the following result. 
\begin{thm}[Corollary~\ref{cor:kr:trc}, Theorem~\ref{lem:tate:frac}]Let $K$ be a nonarchimedean field and let $(r_1,\dots, r_n)\in \mathbb{R}_{>0}^{n}$ be a polyradius. Then \begin{equation*}\dimtr(\Spa(K_{r_{1}, \ldots, r_{n}}, K_{r_1,\dots, r_n}^{\circ}))=\dimtr(\Spa(K_{r_1,\dots, r_{n}}^{\perfd}, K_{r_1,\dots, r_n}^{\perfd\circ}))=n.\end{equation*}\end{thm}
In Section~\ref{sec:proof}, we prove Theorem~\ref{thm:dim:noteq}.
In Section~\ref{sec:open:questions}, we comment on what more refined notions of dimension
might make sense for perfectoid Tate algebras.

\subsection{Acknowledgements}The authors wish to thank their PhD advisor, Kiran Kedlaya, for numerous helpful conversations, advice, and support. They also thank Peter Scholze for a helpful conversation on the results of this paper and Ian Gleason for his comments on our previous work \cite{Dine-Garzella25}. 

This paper was completed when the first-named author was supported by the Jane Street Graduate Research Fellowship 2026. 

\subsection{AI Disclosure}

AI assistants were used to search for, explain, and clarify
various statements from the literature, as well as to transcribe handwritten notes into LaTeX,
perform mechanical LaTeX tasks, and proofread the authors' arguments.
The mathematical ideas are due to the authors, as is the writing.

\section{Preliminaries}

Let $L/K$ be a nonarchimedean field extension. We recall some definitions due to Temkin in \cite{Temkin21} and to Huber in \cite{Huber96}. For any nonarchimedean field $K$, we denote by $\widehat{\overline{K}}$ the completion of an algebraic closure of $K$. 
\begin{defn}[\cite{Temkin21}, \S~2.1.4]
An element $t\in L$ is \emph{topgebraic} over $K$ if
\begin{equation*}K(t)\subseteq \widehat{\overline{K}}\end{equation*}as subfields of $\widehat{\overline{L}}$. Otherwise, $t$ is
\emph{topologically transcendental} over $K$. 

More generally, a valued field extension $L/K$ is called \emph{topgebraic} (and $L$ is said to be topgebraic over $K$) if $L$ is contained in a completed algebraic closure of $K$. 
\end{defn}
\begin{defn}[\cite{Temkin21}, \S~2.1.6]A subset $S\subseteq L$ is \emph{topgebraically independent} if every $x\in S$ is topologically
transcendental over $\widehat{K(S\setminus \{x\})}$, i.e., if
\begin{equation*}
K(S\setminus \{x\})(x)\not\subseteq \widehat{\overline{K(S\setminus \{x\})}},
\end{equation*}regarded as subfields of $\widehat{\overline{L}}$. \end{defn}
\begin{defn}[\emph{Topological transcendence degree}, \cite{Temkin21}, \S~2.1.6]We define the topological transcendence degree of $L/K$ as
\begin{equation*}
\toptrdeg(L/K)=\sup\{\vert S\vert : S\subseteq L \text{ topgebraically
independent over } K\}.
\end{equation*}
\end{defn}
\begin{defn}[\cite{Temkin21}, \S~2.1.8]
A subset $S\subseteq L$ \emph{topgebraically generates} $L$ (or is a topgebraically generating set of $L$) if $L$ is
contained in the completed algebraic closure of $\widehat{K(S)}$, i.e., if
$L$ is topgebraic over $\widehat{K(S)}$.
\end{defn}
\begin{defn}[Topological transcendence basis, \cite{Temkin21}, \S2.1.8]If there exists a subset $S$ of $L$ which is both topgebraically generating and topgebraically independent over $K$, then such a subset is called a topological transcendence basis of $L$ over $K$. \end{defn} 
\begin{defn}[\emph{Topgebraic generating degree}, \cite{Temkin21}, \S~2.1.9]We define the topgebraic generating degree of $L/K$ as
\begin{equation*}
\Toptrdeg(L/K)=\min\{\vert S\vert : S\subseteq L \text{ topgebraically
generating}\}.
\end{equation*}
\end{defn}

\begin{defn}[\cite{Huber96}, Definition~1.8.2] Let $L/K$ be a nonarchimedean field extension. Set \begin{equation*}K(L)=K(L,1)\coloneqq{} \textrm{the topological closure of the algebraic closure of}~K~\textrm{in}~L,\end{equation*}and define inductively: 
\begin{equation*}
K(L,n+1)\coloneqq K(L,n)(L).
\end{equation*}We then define 
\begin{align*}
\trc(L/K)\coloneqq \min\{\trdeg(K^{\prime}/K)\mid{} & K\subseteq K^{\prime}\subseteq L
\text{ such that}\\
& K^{\prime}(L,n)=L \text{ for some } n\}.
\end{align*}
We call this quantity \emph{Huber's topological transcendence degree} of $L/K$. 
\end{defn}
Note that $K(L)\subseteq \widehat{\overline{K}}$, i.e., $K(L)$ is topgebraic over $K$. By definition, we have\begin{equation*}
K(L,n+1)\subseteq \widehat{\overline{K(L,n)}}\subseteq \dots \subseteq
\widehat{\overline{K}},
\end{equation*}
so $K(L,n)$ is topgebraic over $K$ for all $n$.

\begin{lemma}\label{First bound}For any nonarchimedean field extension $L/K$, \begin{equation*}\trc(L/K)\geq \Toptrdeg(L/K)\geq \toptrdeg(L/K).
\end{equation*}
\end{lemma}
\begin{proof}The second inequality is the content of \cite{Temkin21}, Theorem~3.2.1. To prove the first inequality, assume $\trc(L/K)<\infty$ and let \begin{equation*}K\subseteq K^{\prime}\subseteq L\end{equation*}be a subextension such that $K^{\prime}(L,n)=L$ for some $n$ and \begin{equation*}\trdeg(K^{\prime}/K)=\trc(L/K).\end{equation*}
Let $S$ be a transcendence basis of $K'$ over $K$, so $\vert S\vert=\trc(L/K)$. Then $L$ is topgebraic over $\widehat{K(S)}$, so $S$ is a topgebraically generating set of $L$ with \begin{equation*}\vert S\vert\leq \trc(L/K).\end{equation*}It follows that $\Toptrdeg(L/K)\leq \trc(L/K)$, as claimed.
\end{proof}
\begin{cor}\label{First equality}If $\trc(L/K)<\infty$, then $\Toptrdeg(L/K)=\toptrdeg(L/K)$, the nonarchimedean extension $L/K$ possesses a topological transcendence basis and the cardinality of any such basis is $\toptrdeg(L/K)$. \end{cor}
\begin{proof}If $\trc(L/K)<\infty$, then $\Toptrdeg(L/K)<\infty$ by Lemma~\ref{First bound}. The rest follows from this by \cite{Temkin21}, Theorem~3.2.3.  
\end{proof} 
\begin{lemma}\label{Algebraically closed}
If $L$ is algebraically closed, then
\begin{equation*}
\trc(L/K)=\Toptrdeg(L/K).
\end{equation*}
\end{lemma}
\begin{proof}By Lemma~\ref{First bound}, we only need to prove that $\trc(L/K)\leq\Toptrdeg(L/K)$. If $\Toptrdeg(L/K)=\infty$, there is nothing to prove. Suppose that $\Toptrdeg(L/K)<\infty$ and let $S$ be a topgebraically generating set of $L$ with \begin{equation*}\Toptrdeg(L/K)=\vert S\vert.\end{equation*}Since $L$ is algebraically closed, $L$ is the completed algebraic closure of $K^{\prime}=K(S)$. If all elements of $S$ are algebraic over $K$, then $L$ is the completed algebraic closure of $K$ and $\trc(L/K)=0$. Hence we may assume that $S$ is not algebraic over $K$. Choose a maximal algebraically independent subset $S^{\prime}\subseteq S$ and set $K^{\prime\prime}=K(S^{\prime})$. Since $K^{\prime}$ is algebraic over $K^{\prime\prime}=K(S^{\prime})$ and $L$ is the completed algebraic closure of $K^{\prime}$, we see that $L$ is also the completed algebraic closure of $K^{\prime\prime}$. Then \begin{equation*}\trc(L/K)\leq \trdeg(K^{\prime\prime}/K)=\vert S^{\prime}\vert\leq \vert S\vert.\end{equation*}This concludes the proof. \end{proof}

\section{Topological transcendence degree of Berkovich's fields}\label{sec:Berkovich's fields}

We call a tuple of positive real numbers $(r_1,\dots, r_l)\in \mathbb{R}_{>0}^{l}$ a polyradius of length $l$. If $K$ is a fixed nonarchimedean field, we say that a polyradius $(r_1,\dots, r_l)$ is $K$-free (or, if $K$ is understood from the context, simply free) if $r_1,\dots, r_l$ are linearly independent in the ``multiplicative'' $\mathbb{Q}$-vector space \begin{equation*}(\mathbb{R}_{>0}/\sqrt{\vert K^{\times}\vert})\otimes_{\mathbb{Z}}\mathbb{Q},\end{equation*}i.e., if \begin{equation*}r_{1}^{q_{1}}\cdot~\dots~\cdot r_{l}^{q_{l}}\not\in \vert K^{\times}\vert\end{equation*}for all rational numbers $q_1,\dots, q_l$ with $(q_1,\dots, q_l)\neq (0,\dots, 0)$. 
\begin{defn}[Berkovich's fields]For any nonarchimedean field $K$ and any polyradius $(r_1,\dots, r_l)$ ($K$-free or otherwise), we denote by $K_{r_1,\dots, r_l}$ the completion of the rational function field $K(x_1,\dots, x_l)$ with respect to the $(r_1,\dots, r_l)$-Gauss norm, given on polynomials by \begin{equation*}\vert \sum_{i=(i_1,\dots, i_{l})}c_{i}x_{1}^{i_{1}}\cdot~\dots~\cdot x_{l}^{i_{l}}\vert=\max_{i}\vert c_{i}\vert r_{1}^{i_{1}}\cdot~\dots~\cdot r_{l}^{i_{l}}.\end{equation*}We also denote by $K_{r_1,\dots, r_l}^{\perfd}$ the analogous completion of $K(x_{1}^{1/p^{\infty}},\dots, x_{l}^{1/p^{\infty}})$. \end{defn}
If the polyradius $(r_1,\dots, r_l)$ is free, then the maximum in the definition of the Gauss norm is always unique and the fields $K_{r_1,\dots, r_l}$ and $K_{r_1,\dots, r_l}^{\perfd}$ can be described as certain algebras of Laurent series (see, for example, \cite{Dine-Garzella25}, Lemma~2.5). In this case, the notation $K_{r_1,\dots, r_l}$ goes back to Berkovich in \cite{Berkovich}. For a free polyradius, we also have the following statement about residue fields. For any nonarchimedean field $L$ we denote the residue field of $L$ by $\widetilde{L}$. More generally, for any Banach ring $A$, we denote by $\widetilde{A}$ its reduction \begin{equation*}\widetilde{A}=A^{\circ}/A^{\circ\circ},\end{equation*}where $A^{\circ}$ is the subring of power-bounded elements in $A$ and $A^{\circ\circ}$ is the set of topologically nilpotent elements.  
\begin{lemma}\label{Residue fields}Let $K$ be any nonarchimedean field and let $(r_1,\dots, r_n)\in \mathbb{R}_{>0}^{n}$ be a $K$-free polyradius. Then \begin{equation*}\widetilde{K_{r_1,\dots, r_n}}=\widetilde{K}\end{equation*}and\begin{equation*}\widetilde{K_{r_1,\dots, r_n}^{\perfd}}=\widetilde{K}.\end{equation*} \end{lemma}
\begin{proof}Let \begin{equation*}L\in \{K_{r_1,\dots, r_n}, K_{r_1,\dots, r_n}^{\perfd}\}\end{equation*}and consider the canonical map of residue fields $\widetilde{K}\to\widetilde{L}$. We want to prove that this map is surjective. By \cite{Dine-Garzella25}, Lemma~2.5, the nonarchimedean fields $K_{r_1,\dots, r_n}$ and $K_{r_1,\dots, r_n}^{\perfd}$ can be described as the Banach $K$-algebras of Laurent series 
\begin{equation*}K_{r_{1},\dots, r_{n}} = \left\{ \sum_{i=(i_1,\dots, i_{n}) \in \mathbb{Z}^{n}}^{} a_{i}x_{1}^{i_{1}}\cdot~\dots~\cdot x_{n}^{i_{n}} \mid a_{i} \in K, a_{i}r^{i} \xrightarrow{} 0 \text{ as } \vert i\vert
    \xrightarrow{} \infty \right\}\end{equation*}and \begin{equation*}K_{r_1,\dots, r_n}^{\perfd}=\left\{ \sum_{i=(i_1,\dots, i_n) \in \mathbb{Z}[1 / p]^{n}} a_{i}x_{1}^{i_{1}}\cdot~\dots~\cdot x_{n}^{i_{n}} \mid \text{ for every } \epsilon>0, \{i\mid \vert a_{i}\vert r_{1}^{i_{1}}\cdot~\dots~\cdot r_{n}^{i_{n}}\geq\epsilon\} \text{ is finite } \right\}\end{equation*}with the norms being the $(r_1,\dots, r_n)$-Gauss norms. Let \begin{equation*}f=\sum_{i=(i_1,\dots, i_n)} a_{i}x_{1}^{i_{1}}\cdot~\dots~\cdot x_{n}^{i_{n}}\in L\end{equation*}be such a Laurent series of norm \begin{equation*}\vert f\vert=\max_{i=(i_1,\dots, i_n)}\vert a_{i}\vert r_{1}^{i_{1}}\cdot~\dots~\cdot r_{n}^{i_{n}}=1.\end{equation*}Since $(r_1,\dots, r_n)$ is a free polyradius, \begin{equation*}r_{1}^{-i_{1}}\cdot~\dots~r_{n}^{-i_{n}}\not\in \sqrt{\vert K^{\times}\vert}\end{equation*}for any $(i_1,\dots, i_n)\neq (0,\dots, 0)$. It follows that $\vert a_{0}\vert=1$ and \begin{equation*}\vert a_{i}\vert r_{1}^{i_{1}}\cdot~\dots~\cdot r_{n}^{i_{n}}<1\end{equation*}for all $(i_1,\dots, i_n)\neq (0,\dots, 0)$. Consequently, \begin{equation*}f\in a_{0}+L^{\circ\circ}\subseteq K^{\circ}+L^{\circ\circ}.\end{equation*}Since $f\in L$ was arbitrary, this shows that the map $\widetilde{K}=K^{\circ}/K^{\circ\circ}\to \widetilde{L}=L^{\circ}/L^{\circ\circ}$ is surjective, as claimed. \end{proof}
We also record the following consequence of Lemma~\ref{Residue fields} for later use. 
\begin{cor}\label{Reductions}Let $K$ be a nonarchimedean field, let $(r_1,\dots, r_n)\in \mathbb{R}_{>0}^{n}$ be a $K$-free polyradius and let $L\in \{K_{r_1,\dots, r_n}, K_{r_1,\dots, r_n}^{\perfd}\}$. Let \begin{equation*}\varphi: A\to L\end{equation*}be any continuous $K$-algebra homomorphism from a Banach $K$-algebra $A$. Then the induced map of reductions \begin{equation*}\widetilde{\varphi}: \widetilde{A}\to \widetilde{L}\end{equation*}is surjective. In other words, \begin{equation*}\beta\in \varphi(A^{\circ})+L^{\circ\circ}\end{equation*}for all $\beta\in L$ with $\vert \beta\vert=1$. \end{cor}
\begin{proof}The map $\widetilde{\varphi}$ is a map of $\widetilde{K}$-algebras, so the image of $\widetilde{\varphi}$ is a $\widetilde{K}$-subalgebra of $\widetilde{L}$. But $\widetilde{L}=\widetilde{K}$ by Lemma~\ref{Residue fields}. \end{proof}
By contrast, if the polyradius $(r_1,\dots, r_n)$ is not free, then the residue field of $K_{r_1,\dots, r_n}$ has positive transcendence degree over the residue field of $K$. \begin{lemma}\label{Rational}Let $K$ be a nonarchimedean field, let $r\in\sqrt{\vert K^{\times}\vert}$ and let $c\in K^{\times}$, $k\geq1$ such that $r^{k}=\vert c\vert$. Let $x\in K_{r}$ be an element such that $K_{r}$ is the completion of $K(x)$ for the $r$-Gauss norm. Then the image \begin{equation*}\widetilde{y}\in \widetilde{K_{r}}\end{equation*}of $y=\frac{x^{k}}{c}$ in the residue field of $K_{r}$ is transcendental over the residue field $\widetilde{K}$ of $K$. \end{lemma}
\begin{proof}Suppose that \begin{equation*}\sum_{j}\widetilde{a_{j}}\widetilde{y}^{j}=0, \widetilde{a_{j}}\in \widetilde{K}\setminus\{0\},\end{equation*}is an algebraic relation over $\widetilde{K}$, and lift each $\widetilde{a_{j}}$ to an element $a_{j}\in K$ with $\vert a_{j}\vert=1$. 
Then \(\sum_{j}^{} a_{j}y^{j} \) is in the maximal ideal \(K_{r}^{\circ\circ}\) of \(K_{r}^{\circ}\),
thus \begin{equation*}\vert \sum_{j}a_{j}y^{j}\vert_{r}<1.\end{equation*}But, by definition of the $r$-Gauss norm, \begin{equation*}\vert \sum_{j}a_{j}y^{j}\vert_{r}=\vert \sum_{j}a_{j}c^{-j}x^{jk}\vert_{r}=\max_{j}\vert a_{j}\vert \vert c\vert^{-j}r^{jk}=1,\end{equation*}which produces the desired contradiction. \end{proof}
\begin{lemma}\label{One-variable case}Let $K$ be a nonarchimedean field and let $r>0$. Let $K_{r}$ be the completion of the purely transcendental extension field $K(x)$ for the $r$-Gauss norm. Then $x$ is topologically transcendental over $K$ and $\toptrdeg(K_r/K)=1$. \end{lemma}
\begin{proof}The element $x$ clearly generates $K_r$ topgebraically, so it suffices to prove that it is topologically transcendental over $K$. In other words, we need to prove that $K(x)$ is not contained in a completed algebraic closure of $K$. 
If $K(x)$ were contained in a completed algebraic closure, then the value group $\vert K_{r}^{\times}\vert$ would be contained in the divisible closure of $\vert K^{\times}\vert$ and its residue field would be algebraic over that of $K$. 
If \begin{equation*}r\not\in \sqrt{\vert K^{\times}\vert},\end{equation*}then \begin{equation*}\vert K_{r}^{\times}\vert=\vert K^{\times}\vert\cdot r^{\mathbb{Z}}\not\subseteq \sqrt{\vert K^{\times}\vert}.\end{equation*}
Otherwise, $r^{k}=\vert c\vert$ for some $c\in K^{\times}$ and $k\geq1$. But then the image of $y=\frac{x^{k}}{c}$ in the residue field of $K_{r}$ is transcendental over the residue field of $K$, by Lemma~\ref{Rational}.\end{proof}
\begin{thm}\label{Main theorem}
For any nonarchimedean field $K$ and any polyradius $(r_1,\dots, r_l)\in \mathbb{R}_{>0}^{l}$, 
\begin{equation*}
\toptrdeg(K_{r_1,\dots, r_l}/K)=\Toptrdeg(K_{r_1,\dots, r_l}/K)=l.
\end{equation*}Moreover, any algebraically independent subset $S=\{x_1,\dots, x_l\}$ such that $K_{r_1,\dots, r_l}$ is the completion of $K(x_1,\dots, x_l)$ for the $(r_1,\dots, r_l)$-Gauss norm is a topological transcendence basis of $K_{r_1,\dots, r_l}$. 
\end{thm}
\begin{proof}The second assertion clearly implies the first. It is clear that the subset $S=\{x_1,\dots, x_l\}$ of $K_{r_1,\dots, r_l}$ is topgebraically generating and we only need to show that this subset is also topgebraically independent. That is, we need to show that for every index $i\in\{1,\dots, l\}$, the field $K_{r_1,\dots, r_l}$ is not contained in the completed algebraic closure of $L=K_{r_1,\dots, \hat{r_{i}},\dots, r_l}$. We can write the nonarchimedean field $K_{r_1,\dots, r_l}$ as \begin{equation*}K_{r_1,\dots, r_l}=L_{r_{i}},\end{equation*}so, replacing $K$ with $L$, it suffices to prove the assertion in the case $l=1$. This case was already treated in Lemma~\ref{One-variable case}.  \end{proof}

\begin{cor}
    \label{cor:kr:trc}
For any nonarchimedean field $K$ and any polyradius $(r_1,\dots, r_l)\in \mathbb{R}_{>0}^{l}$,
\begin{equation*}
\trc(K_{r_1,\dots, r_l}/K)=l.
\end{equation*}
\end{cor}

\begin{proof}
By Lemma~\ref{First bound}, we have the inequality 
\begin{equation*}
\trc(K_{r_1,\dots, r_l}/K)\geq \Toptrdeg(K_{r_1,\dots, r_l}/K).
\end{equation*}
By Theorem~\ref{Main theorem}, this implies that $\trc(K_{r_1,\dots, r_l}/K)\geq l$. On the other hand, for any algebraically independent subset $S=\{x_1,\dots, x_l\}$ such that $K_{r_1,\dots, r_l}$ is the completion of $K(x_1,\dots, x_l)$ for the $(r_1,\dots, r_l)$-Gauss norm, we have \begin{equation*}
K(x_1,\dots, x_l)(K_{r_1,\dots, r_l}, 1)=K_{r_1,\dots, r_l},
\end{equation*}
so $\trc(K_{r_1,\dots, r_l}/K)\leq \trdeg(K(x_1,\dots, x_l)/K)=l$.
\end{proof}
Similarly, we get the following result about the perfectoid fields $K_{r_1,\dots, r_n}^{\perfd}$. 
\begin{thm}
    \label{lem:tate:frac}
    Let \(K\) be a nonarchimedean field, let $n\geq1$, $m\geq0$ and let $(r_1,\dots, r_n)\in \mathbb{R}_{>0}^{n}$ be a polyradius. 
    Let \(L\) be the completion of the fraction field of \(K_{r_{1}, \ldots, r_{n}}^{\perfd}\langle
    T_{1}^{1/p^{\infty}}, \ldots, T_{m}^{1/p^{\infty}} \rangle \).
    Then, \(\trc(L / K) = n+m.\)
\end{thm}
\begin{proof}
    Let $\{x_1,\dots, x_n\}$ be a finite algebraically independent subset of $K_{r_1,\dots, r_n}$ such that $K_{r_1,\dots, r_n}$ is the completion of $K(x_1,\dots, x_n)$ for the $(r_1,\dots, r_n)$-Gauss norm. Let \(K^{\prime} = K(x_{1}, \ldots, x_{n}, T_{1}, \ldots, T_{m})\). The absolute value on the completion of $K^{\prime}$ extending the Gauss norm on $K_{r_1,\dots, r_n}\langle T_1,\dots, T_m\rangle$ is the $(r_1,\dots, r_n, 1, \dots, 1)$-Gauss norm.  
    By Theorem~\ref{Main theorem}, the subset $\{x_1,\dots, x_n, T_1,\dots, T_m\}$ is a topological transcendence basis of (the completion of) $K^{\prime}$. 
    In particular, this is a
    topgebraically independent subset of \(L\).
    Thus, \(\trc(L / K)\) is at least \(n + m\)
    by Lemma~\ref{First bound}.
    
    On the other hand, \(K^{\prime}\) has transcendence degree \(n+m\),
    and we will show that \(L = K^{\prime}(L,3)\), or
    equivalently \(L \ins K^{\prime}(L,3)\).
    First, 
    \(K(x_{1}^{1 / p^{\infty}}, \ldots, x_{n}^{1 / p^{\infty}}, T_{1}^{1 / p^{\infty}}, \ldots,
    T_{m}^{1 / p^{\infty}}) \ins K^{\prime}(L,1)\), 
    as it is inside the algebraic closure of \(K^{\prime}\) in \(L\).
    Furthermore, 
    \(K_{r_{1}, \ldots, r_{n}}^{\perfd}(T_{1}^{1 / p^{\infty}}, \ldots, T_{m}^{1 / p^{\infty}}) \ins
    K^{\prime}(L,1)\),
    as that lives in the topological closure of 
    \(K(x_{1}^{1 / p^{\infty}}, \ldots, x_{n}^{1 / p^{\infty}}, T_{1}^{1 / p^{\infty}}, \ldots,
    T_{m}^{1 / p^{\infty}})\).
    Next, 
    we see that 
    \(K_{r_{1}, \ldots, r_{n}}^{\perfd}\langle T_{1}^{1 / p^{\infty}}, \ldots, T_{m}^{1 / p^{\infty}}\rangle
    \ins K^{\prime}(L,2)\), 
    and thus its fraction field also lies in \(K^{\prime}(L,2)\).
    Thus, we see that the topological closure of this fraction
    field, which is \(L\), is a subset of \(K^{\prime}(L,3)\).
\end{proof}



\section{Proof of Theorem \ref{thm:dim:noteq}}
\label{sec:proof}

Let $X\to \Spa(K, K^{\circ})$ be an adic space over a nonarchimedean field $K$. Following \cite{Huber96}, Definition~1.8.4, we define $\dimtr(X)$ to be the supremum over $\trc(\kappa(x)/K)$ where $x$ ranges over the points $x\in X$ and where $\kappa(x)$ denotes the completed residue field of $X$ at the point $x$. 
\begin{defn}[Iterated transcendence degree]We call $\dimtr(X)$ the iterated transcendence degree of $X$ over $K$. \end{defn}
We also denote by $\dim(X)$ the usual Krull dimension of the underlying locally spectral space of the adic space $X$, i.e., the supremum over the lengths of chains of specializations in $X$. We use the following somewhat stronger form of one of our main theorems from \cite{Dine-Garzella25}. 
\begin{thm}[Stronger form of \cite{Dine-Garzella25}, Theorem~6.1]\label{Old theorem}Let $K$ be a perfectoid field of characteristic $p$. Then for every integer $n\geq1$ and every $K$-free polyradius $(r_1,\dots, r_n)\in (0, 1]^{n}$, there exists a surjective map \begin{equation*}T_{n+2, K}^{\perfd\circ}\to K_{r_1,\dots, r_n}^{\perfd\circ}.\end{equation*}\end{thm}
\begin{proof}We deduce this somewhat stronger statement from the results we proved in \cite{Dine-Garzella25}. Consider the map \begin{equation*}\varphi: R=T_{n+2, K}^{\perfd}\to K_{r_1,\dots, r_n}^{\perfd}\end{equation*}from \cite{Dine-Garzella25}, Corollary~6.19. We want to prove that this map restricts to a surjective map on subrings of power-bounded elements. 

Let $\beta\in K_{r_1,\dots, r_n}^{\perfd\circ}$. We want to prove that $\beta\in \varphi(R^{\circ})$. If $\vert \beta\vert=1$, then Corollary~\ref{Reductions} tells us that \begin{equation*}\beta\in \varphi(R^{\circ})+K_{r_1,\dots, r_n}^{\perfd\circ\circ}.\end{equation*}Therefore, we may assume that $\vert \beta\vert<1$. As we explained in the proof of \cite{Dine-Garzella25}, Corollary~6.19, the map $\varphi$ satisfies the assumptions of loc.~cit., Proposition~6.14. These assumptions imply that for some $s\in (0, 1)$, the image $\varphi(R^{\circ})$ of $R^{\circ}$ contains a $(q, s)$-adapted element for every $q\in \mathbb{Z}[1/p]^{n}$ in the sense of loc.~cit., Definition~6.2. Since $\vert \beta\vert<1$, we can choose some $m\geq0$ such that $\vert \beta^{p^{m}}\vert\leq s$. Then the proof of loc.~cit., Proposition~6.6 (which, in its original form, is due to Gleason, see \cite{gleason-2022-perfectoid-Nullstellensatz}, proof of Lemma~4.3) produces an element $f\in R^{\circ}$ with $\varphi(f)=\beta^{p^{m}}$. But then $\varphi(f^{1/p^{m}})=\beta$. \end{proof}
\begin{thm}\label{Main theorem 2}Let $K$ be a perfectoid field of characteristic $p$ with countable value group. For every integer $n\geq1$, there exists a Zariski-closed subset $Z$ of the $(n+2)$-dimensional perfectoid closed unit ball $\mathbb{B}_{K}^{n+2,\perfd}$ over $K$ such that $\dim(Z)=0$ and $\dimtr(Z)=n$.\end{thm}
\begin{proof}
Since the value group of $K$ is countable, so is the divisible closure of the value group. It follows that the $\mathbb{Q}$-vector space \begin{equation*}(\mathbb{R}_{>0}/\sqrt{\vert K^{\times}\vert})\otimes_{\mathbb{Z}}\mathbb{Q} \end{equation*}is of uncountable dimension. In particular, \begin{equation*}\dim_{\mathbb{Q}}(\mathbb{R}_{>0}/\sqrt{\vert K^{\times}\vert})\geq n.\end{equation*}This means that there exists a $K$-free polyradius $(r_1,\dots, r_n)\in \mathbb{R}_{>0}^{n}$ of length $n$. Up to multiplying $r_1,\dots, r_n$ by a sufficiently small element in $\vert K^{\times}\vert$, we may assume that $r_1,\dots, r_n\in (0, 1]$. In what follows, we fix such a free polyradius $(r_1,\dots, r_n)\in (0, 1]^{n}$. 

Set $L=K_{r_1,\dots, r_n}^{\perfd}$. By Theorem~\ref{Old theorem}, there exists a surjective map \begin{equation*}T_{n+2, K}^{\perfd\circ}\to L^{\circ},\end{equation*}so $\Spa(L, L^{\circ})$ is a Zariski-closed subset of the $(n+2)$-dimensional perfectoid closed unit ball $\mathbb{B}_{K}^{n+2,\perfd}$. By the special case $m=0$ of Theorem~\ref{lem:tate:frac}, \begin{equation*}\trc(L/K)=n.\end{equation*}But $\Spa(L, L^{\circ})$ is clearly of Krull dimension $0$. 
\end{proof}
\begin{lemma}\label{Trivial upper bound}Let $n\geq1$ be an integer, let $K$ be a perfectoid field and let $Z$ be a Zariski-closed subset of $\mathbb{B}_{K}^{n,\perfd}$. Then for every $x\in Z$, the completed residue field $\kappa(x)$ at $x$ satisfies $\trc(\kappa(x)/K)\leq n$. In particular, $\dimtr \mathbb{B}_{K}^{n,\perfd}\leq n$. \end{lemma}
\begin{proof}We can view $x$ as a point of $\mathbb{B}_{K}^{n,\perfd}$. Recall that $\kappa(x)$ is a completion of the fraction field of $T_{n, K}^{\perfd}/\ker(x)$. Let $f_1,\dots, f_n$ be the images of the variables $X_1, \dots, X_n$ in $T_{n, K}^{\perfd}/\ker(x)$. 
Let $S$ be a maximal algebraically independent subset of $\{f_1,\dots, f_n\}$ (could possibly be empty). Then the topological closure of the algebraic closure of $K(S)$ in $\kappa(x)$ contains the image $T_{n, K}^{\perfd}/\ker(x)$ of $T_{n, K}^{\perfd}$ in $\kappa(x)$. 
Hence it also contains the field of fractions of $T_{n, K}^{\perfd}/\ker(x)$. It follows that $\kappa(x)$ is equal to the topological closure of the algebraic closure of $K(S)$ in $\kappa(x)$. 
Thus \begin{equation*}\kappa(x)=K(S)(\kappa(x), 1).\end{equation*}In particular, $\trc(\kappa(x)/K)\leq \vert S\vert\leq n$. \end{proof}
\begin{cor}\label{Dimension of the perfectoid closed unit ball}Let $K$ be a perfectoid field. Then \begin{equation*}\dim\mathbb{B}_{K}^{n,\perfd}=\dimtr\mathbb{B}_{K}^{n,\perfd}=n\end{equation*}for all integers $n\geq1$. \end{cor}
\begin{proof}The inequality $\dim X\leq \dimtr X$ holds for any adic space $X$ over $\Spa(K, K^{\circ})$ by \cite{Huber96}, Lemma~1.8.5i) (which is a consequence of Abhyankar's inequality). Combining this with Lemma~\ref{Trivial upper bound}, we get \begin{equation*}\dim\mathbb{B}_{K}^{n,\perfd}\leq\dimtr\mathbb{B}_{K}^{n,\perfd}\leq n.\end{equation*}On the other hand, there is a surjective, continuous and closed specialization map \begin{equation*}\mathbb{B}_{K}^{n,\perfd}=\Spa(T_{n, K}^{\perfd}, T_{n, K}^{\perfd\circ})\to \Spf(T_{n, K}^{\perfd\circ})=\Spec(\widetilde{K}[X_{1}^{1/p^{\infty}}, \dots, X_{n}^{1/p^{\infty}}]),\end{equation*}so $\dim\mathbb{B}_{K}^{n,\perfd}\geq n$. \end{proof} 
\begin{thm}
\label{Main theorem 3}
Let $K$ be a perfectoid field of characteristic $p$ whose value group is countable. 
For every pair of integers $n\geq1$, $m\geq0$, there exists a Zariski-closed subset $Z$ of the $(n+m+2)$-dimensional
perfectoid closed unit ball $\mathbb{B}_{K}^{n+m+2,\perfd}$ over $K$ such that $\dim(Z)=m$ and
$\dimtr(Z)\geq n+m$.
\end{thm}
\begin{proof}
By the proof of Theorem~\ref{Main theorem 2}, there exists a $K$-free polyradius $(r_1,\dots, r_n)\in \mathbb{R}_{>0}^{n}$ such that \(\Spa(K_{r_{1}, \ldots, r_{n}}^{\perfd}, K_{r_{1},
    \ldots, r_{n}}^{\perfd\circ})\) is a Zariski-closed subset of \(\mathbb{B}_{K}^{n+2,\perfd}\).
    Thus, \begin{equation*}Z = \Spa(K_{r_{1}, \ldots, r_{n}}^{\perfd}\langle T_{1}^{1 / p^{\infty}}, \ldots, T_{m}^{1 /
        p^{\infty}}\rangle, K_{r_{1},
    \ldots, r_{n}}^{\perfd\circ}\langle T_{1}^{1 / p^{\infty}}, \ldots, T_{m}^{1 / p^{\infty}}\rangle)\end{equation*}is a Zariski-closed subset of $\mathbb{B}_{K}^{n+m+2,\perfd}$. 
    Then \(L\) as in Theorem~\ref{lem:tate:frac} is a residue field of \(Z\) with
    \(\trc(L / K) = n+m\).
\end{proof}


\section{Open questions}
\label{sec:open:questions}

In Berkovich's theory, which allows for non-strictly affinoid algebras
(i.e., quotients of \(K\langle X_{1} / r_{1}, \ldots, X_{n} / r_{n}\rangle \) 
for \(r_{i}\) possibly not in the divisible closure of the value group of \(K\)),
the notion of dimension is slightly more refined. Before we state the definition, we first record the proofs of two classical facts used in this context. 
\begin{lemma}\label{Common field extensions}Let $L_{1}/K$ and $L_{2}/K$ be two nonarchimedean field extensions. Then there exists a common nonarchimedean field extension $L$ of both $L_{1}$ and $L_{2}$. \end{lemma}
\begin{proof}Consider the Banach ring $L_{1}\widehat{\otimes}_{K}L_{2}$. By \cite{Gruson}, Théorème~1~\no4, the map \begin{equation*}L_{1}\otimes_{K}L_{2}\to L_{1}\widehat{\otimes}_{K}L_{2}\end{equation*}is injective, so, in particular, \begin{equation*}L_{1}\widehat{\otimes}_{K}L_{2}\neq 0.\end{equation*}By \cite{Berkovich}, Theorem~1.2.1, there exists a bounded multiplicative seminorm $v$ on $L_{1}\widehat{\otimes}_{K}L_{2}$. Consider the integral domain $A=(L_{1}\widehat{\otimes}_{K}L_{2})/\ker(v)$ equipped with the absolute value induced by $v$, which we again denote by $v$, and extend the absolute value to the field of fractions of $A$. Let $L$ be the completion of $\Frac(A)$ with respect to this absolute value. Since $v$ was bounded with respect to the norm on $L_{1}\widehat{\otimes}_{K}L_{2}$, we have bounded ring maps \begin{equation*}L_{i}\to A\hookrightarrow L\end{equation*}for $i=1, 2$ which are injective since their sources are fields. Thus the absolute value $v$ on $L$ restricts to absolute values on $L_{1}$, $L_{2}$ which are bounded with respect to their given absolute values. Since the Berkovich spectrum of any nonarchimedean field consists of a single point (\cite{Kedlaya18}, Remark~2.11), we conclude that $v$ restricts to the given absolute values on $L_{1}$ and $L_{2}$, i.e., $L$ is a nonarchimedean field extension of both $L_{1}$ and $L_{2}$.  \end{proof}
The following statement is used a lot in proofs in the theory of Berkovich spaces. 
\begin{lemma}[cf.~\cite{Berkovich}, \S2.1]\label{Strictly affinoid}Let $A$ be a $K$-affinoid algebra in the sense of Berkovich. Then there exists a polyradius $(r_1,\dots, r_n)\in \mathbb{R}_{>0}^{n}$ such that $A\widehat{\otimes}_{K}K_{r_1,\dots, r_n}$ is a strictly affinoid $K_{r_1,\dots, r_n}$-algebra.\end{lemma}
\begin{proof}By definition of Berkovich-affinoid algebras, there exist an integer $n\geq0$, a polyradius $(r_1,\dots, r_n)\in \mathbb{R}_{>0}^{n}$ and a strict surjective ring map \begin{equation*}K\langle r_{1}^{-1}X_1,\dots, r_{n}^{-1}X_n\rangle \to A.\end{equation*}By \cite{BGR}, Proposition~2.1.8/6, the map \begin{equation*}K_{r_1,\dots, r_n}\langle r_{1}^{-1}X_{1},\dots, r_{n}^{-1}X_{n}\rangle=K\langle r_{1}^{-1}X_1,\dots, r_{n}^{-1}X_n\rangle\widehat{\otimes}_{K}K_{r_1,\dots, r_n}\to A\widehat{\otimes}_{K}K_{r_1,\dots, r_n}\end{equation*}is also strict surjective. But the algebra $K_{r_1,\dots, r_n}\langle r_{1}^{-1}X_{1},\dots, r_{n}^{-1}X_{n}\rangle$ is strictly $K_{r_1,\dots, r_n}$-affinoid (by \cite{BGR}, Theorem~6.1.5/4) since all the radii $r_1,\dots, r_n$ belong to the value group of its ground field $K_{r_1,\dots, r_n}$. Hence $A\widehat{\otimes}_{K}K_{r_1,\dots, r_n}$ is a quotient of a strictly $K_{r_1,\dots, r_n}$-affinoid algebra, which makes it strictly $K_{r_1,\dots, r_n}$-affinoid. \end{proof}
One can prove a more precise statement, which extends both the previous lemma and \cite{Berkovich}, Corollary~2.1.6. 
\begin{lemma}\label{Strictly affinoid 2}Let $A$ be a Berkovich-affinoid $K$-algebra and let $L/K$ be a nonarchimedean field extension. If $A\widehat{\otimes}_{K}L$ is strictly $L$-affinoid, then \begin{equation*}\vert A\vert_{\spc}\subseteq \sqrt{\vert L^{\times}\vert}\cup\{0\},\end{equation*}where $\vert\cdot\vert_{\spc}$ is the spectral seminorm on $A$. Conversely, if $A$ is reduced and $\vert A\vert_{\spc}\subseteq \sqrt{\vert L^{\times}\vert}\cup\{0\}$, then $A\widehat{\otimes}_{K}L$ is strictly $L$-affinoid.\end{lemma}
\begin{proof}By \cite{Gruson}, Théorème~1~\no1, applied to the strict injective map $K\hookrightarrow L$ and the complete base extension functor $-\widehat{\otimes}_{K}A$, the map of Banach $K$-algebras \begin{equation*}A\hookrightarrow A\widehat{\otimes}_{K}L\end{equation*}is injective and strict. That is, the map is injective and the norm on $A$ is equivalent to the restriction of the norm on $A\widehat{\otimes}_{K}L$. By \cite{BGR}, Proposition~1.3.1/2, this implies that the spectral seminorm on $A\widehat{\otimes}_{K}L$ restricts to the spectral seminorm on $A$. We conclude that if $A\widehat{\otimes}_{K}L$ is strictly affinoid, then \begin{equation*}\vert A\vert_{\spc}\subseteq \vert A\widehat{\otimes}_{K}L\vert_{\spc}\subseteq \sqrt{\vert L^{\times}\vert}\cup\{0\}\end{equation*}by \cite{BGR}, Proposition~6.2.1/4(ii). 

For the converse implication, we imitate the proof of \cite{Berkovich}, Corollary~2.1.6. Suppose that $A$ is reduced. Let \begin{equation*}\varphi: K\langle r_{1}^{-1}X_{1},\dots, r_{n}^{-1}X_{n}\rangle\to A\end{equation*}be a strict surjection (for some polyradius $(r_1,\dots, r_n)\in\mathbb{R}_{>0}^{n}$). For $i=1,\dots, n$, let \begin{equation*}f_{i}=\varphi(X_{i}), r_{i}^{\prime}=\vert f_{i}\vert_{\spc}.\end{equation*}We may assume that $f_{i}\neq0$ for all $i$. Since $A$ is reduced, the spectral seminorm $\vert\cdot\vert_{\spc}$ on $A$ is a norm by \cite{Berkovich}, Proposition~2.1.4, so $r_{i}^{\prime}\neq 0$ for all $i$. The map \begin{equation*}K[X_1,\dots, X_n]\to K[f_1,\dots, f_n]\subseteq A\end{equation*}given by restriction of $\varphi$ is bounded with respect to the $(r_{1}^{\prime},\dots, r_{n}^{\prime})$-Gauss norm on $K[X_1,\dots, X_n]$ since \begin{equation*}\vert \sum_{i=(i_1,\dots, i_n)}a_{i}f_{1}^{i_{1}}\cdot~\dots~\cdot f_{n}^{i_{n}}\vert_{\spc}\leq \max_{i}\vert a_{i}\vert \vert f_{1}^{i_{1}}\vert_{\spc}\cdot~\dots~\cdot \vert f_{n}^{i_{n}}\vert_{\spc}\leq\max_{i}\vert a_{i}\vert r_{1}^{\prime i_{1}}\cdot~\dots~\cdot r_{n}^{\prime i_{n}}.\end{equation*}On the other hand, since $\varphi$ is bounded, $r_{i}^{\prime}\leq r_{i}$ for all $i=1,\dots, n$ by virtue of \cite{BGR}, Proposition~1.3.1/2. Thus the $(r_{1}^{\prime},\dots, r_{n}^{\prime})$-Gauss norm on the polynomial ring $K[X_1,\dots, X_n]$ is bounded above by the $(r_1,\dots, r_n)$-Gauss norm. It follows by \cite{BGR}, Corollary~1.1.9/6, that the map $K[X_1,\dots, X_n]\to K[f_1,\dots, f_n]$ extends to a strict surjective map \begin{equation*}K\langle r_{1}^{\prime-1}X_{1},\dots, r_{n}^{\prime-1}X_{n}\rangle\to A\end{equation*}(we are using here that, since $A$ is reduced, the topology on $A$ is the topology defined by the spectral norm $\vert\cdot\vert_{\spc}$, by \cite{Berkovich}, Proposition~2.1.4). By \cite{BGR}, Proposition~2.1.8/6, we obtain a strict surjective map \begin{equation*}L\langle r_{1}^{\prime-1}X_{1},\dots, r_{n}^{\prime-1}X_{n}\rangle\to A\widehat{\otimes}_{K}L.\end{equation*}But if \begin{equation*}\vert A\vert_{\spc}\subseteq \sqrt{\vert L^{\times}\vert}\cup\{0\},\end{equation*}then \begin{equation*}r_{1}^{\prime},\dots, r_{n}^{\prime}\in \sqrt{\vert L^{\times}\vert}\end{equation*}and it follows from \cite{BGR}, Theorem~6.1.5/4, that the affinoid $L$-algebra $A\widehat{\otimes}_{K}L$ is strictly $L$-affinoid. \end{proof}
\begin{lemma}\label{Finite maps}Let $R$ be a complete Tate ring and let $A\to A^{\prime}$, $A\to B$ be continuous maps of Banach $R$-algebras (i.e., complete Tate rings receiving a continuous map from $R$) and suppose that $A\to A^{\prime}$ is a finite ring map. Then the induced map \begin{equation*}B\to B^{\prime}=A^{\prime}\widehat{\otimes}_{A}B\end{equation*}is again finite. \end{lemma}
\begin{proof}Indeed, any surjection of $A$-modules \begin{equation*}A^{m}\to A^{\prime}\end{equation*}is strict by Banach's open mapping theorem, so applying the complete base extension functor $-\widehat{\otimes}_{A}B$ and \cite{BGR}, Proposition~2.1.8/6, we obtain a (strict) surjection \begin{equation*}B^{m}\to B^{\prime},\end{equation*}proving that $B^{\prime}$ is finitely generated as a $B$-module. \end{proof}
\begin{lemma}\label{Finite maps 2}Let $K$ be a nonarchimedean field and let $A\to A^{\prime}$ be a finite, injective and strict map of affinoid $K$-algebras in the sense of Berkovich. For any nonarchimedean field extension $L/K$, the map $A\widehat{\otimes}_{K}L\to A^{\prime}\widehat{\otimes}_{K}L$ is again finite, injective and strict.\end{lemma}
\begin{proof}The map in question is finite by Lemma~\ref{Finite maps}. It is injective and strict by \cite{Gruson}, Théorème~1~\no1. \end{proof} 
\begin{lemma}\label{Noether normalization}Let $K$ be a nonarchimedean field, let $A$ be an affinoid $K$-algebra in the sense of Berkovich and let $L_{1}/K$, $L_{2}/K$ be two nonarchimedean field extensions such that $A\widehat{\otimes}_{K}L_{1}$ and $A\widehat{\otimes}_{K}L_{2}$ are strictly affinoid. Then the Krull dimensions of $A\widehat{\otimes}_{K}L_{1}$ and $A\widehat{\otimes}_{K}L_{2}$ are the same. \end{lemma}
\begin{proof}By Noether normalization for strictly affinoid algebras (\cite{BGR}, Corollary~6.1.2/2) there exist integers $d_{i}\geq0$ and finite injective strict homomorphisms \begin{equation*}L_{i}\langle T_1,\dots, T_{d_{i}}\rangle \hookrightarrow A\widehat{\otimes}_{K}L_{i}\end{equation*}for $i=1, 2$. In particular, the Krull dimension of $A\widehat{\otimes}_{K}L_{i}$ (for $i=1, 2$) is $d_{i}$. By Lemma~\ref{Common field extensions}, there exists a common nonarchimedean field extension $L$ of $L_{1}$ and $L_{2}$. Base-changing to $L$, the maps \begin{equation*}L\langle T_1,\dots, T_{d_{i}}\rangle\hookrightarrow A\widehat{\otimes}_{K}L=(A\widehat{\otimes}_{K}L_{i})\widehat{\otimes}_{L_{i}}L\end{equation*}(for $i=1, 2$) are again finite and injective by Lemma~\ref{Finite maps 2}. It follows that $d_{1}=d_{2}$, as desired.\end{proof}
\begin{defn}
    [Berkovich \cite{Berkovich}, paragraph after Proposition~2.3.4, cf.~\cite{Ducros07}, Définition~1.6]
    Let $K$ be a nonarchimedean field and let \(A\) be a $K$-affinoid algebra in the sense of Berkovich. 
    Let $L/K$ be a nonarchimedean field extension such that \begin{equation*}\vert A\vert_{\spc}\subseteq \sqrt{\vert L^{\times}\vert}\cup\{0\}\end{equation*}(such $L/K$ exist by Lemma~\ref{Strictly affinoid} and Lemma~\ref{Strictly affinoid 2}). Then the $K$-analytic dimension \begin{equation*}\dim_{K}(A)\end{equation*}of $A$ is the Krull dimension of $A\widehat{\otimes}_{K}L$. \end{defn}
\begin{lemma}\label{Lemma on analytic dimension}The above definition does not depend on the choice of a nonarchimedean field extension $L/K$ with $\vert A\vert_{\spc}\subseteq \sqrt{\vert L^{\times}\vert}\cup\{0\}$.\end{lemma}
\begin{proof}For every nonarchimedean field extension $L/K$, the quotient map $A\to A/\Nil(A)$ induces a surjective map \begin{equation*}A\widehat{\otimes}_{K}L\to (A/\Nil(A))\widehat{\otimes}_{K}L\end{equation*}by \cite{BGR}, Proposition~2.1.8/6, and the kernel of this map is clearly (locally) nilpotent. Hence $A\widehat{\otimes}_{K}L$ has the same Krull dimension as $(A/\Nil(A))\widehat{\otimes}_{K}L$, and it suffices to prove the assertion of the lemma when $A$ is reduced. In this case, the assertion follows from Lemma~\ref{Strictly affinoid 2} and Lemma~\ref{Noether normalization}. \end{proof}
\begin{example}Let $K$ be any nonarchimedean field and let $(r_1,\dots, r_n)\in \mathbb{R}_{>0}^{n}$ be a $K$-free polyradius of some length $n\geq1$. Then the nonarchimedean field $K_{r_1,\dots, r_n}$ is a Berkovich $K$-affinoid algebra whose $K$-analytic dimension is \begin{equation*}\dim_{K}(K_{r_1,\dots, r_n})=n.\end{equation*}Indeed, since the polyradius is free, $K_{r_1,\dots, r_n}$ is the $K$-affinoid algebra of Laurent series \begin{equation*}K_{r_1,\dots, r_n}=K\langle r_{1}^{-1}X_{1},\dots, r_{n}^{-1}X_{n}, r_{1}X_{1}^{-1},\dots, r_{n}X_{n}^{-1}\rangle\end{equation*}(for a proof of this classical fact, see, for example, \cite{Dine-Garzella25}, Lemma~2.5) and then \begin{equation*}K_{r_1,\dots, r_n}\widehat{\otimes}_{K}K_{r_1,\dots, r_n}\end{equation*}is the analogous Laurent series algebra with coefficients in $K_{r_1,\dots, r_n}$, but now with the radii in the value group of the new ground field $K_{r_1,\dots, r_n}$. In particular, $K_{r_1,\dots, r_n}\widehat{\otimes}_{K}K_{r_1,\dots, r_n}$ is a strictly $K_{r_1,\dots, r_n}$-affinoid algebra of Krull dimension $n$.\end{example}
In the context of perfectoid spaces, 
one may wonder if instead of using the (adic space) Krull dimension naively,
we may take a notion similar to Berkovich's $K$-analytic dimension to be the dimension in the case of Zariski-closed subsets of a perfectoid ball. Thus, we may ask:

\begin{quest}
    Let \(K\) be a perfectoid field and let \begin{equation*}Z\subseteq \mathbb{B}_{K}^{n,\perfd}\end{equation*}be a Zariski-closed subset. Let $L/K$ be an extension of perfectoid fields such that $\vert L^{\times}\vert=\mathbb{R}_{>0}$. 
    Is the Krull dimension of the base change $Z\times_{K}L$ necessarily equal to \(\dimtr(Z)\)?
\end{quest}

We may ask even a weaker question:

\begin{quest}
    Let \(K\) be a perfectoid field and let \begin{equation*}Z\subseteq \mathbb{B}_{K}^{n,\perfd}\end{equation*}be a Zariski-closed subset. Does there exist an extension of perfectoid fields $L/K$ such that the Krull dimension of the base change $Z\times_{K}L$ equals $\dimtr(Z)$? 
\end{quest}

\bibliographystyle{amsalpha}
\bibliography{main.bib}

\end{document}